\documentclass[a4paper,12pt,draft]{article}

\makeatletter

\@setfontsize\normalsize\@xiipt{15.5}%

\renewcommand*{\author}[1]{\gdef\@author{#1}\gdef\@pauthor{{\def\and{ --- }#1}}}
\renewcommand*{\title}[1]{\gdef\@title{#1}\gdef\@ptitle{#1}}

\if@twoside         
  \def\ps@draft{%
    \def\@oddfoot{\small\null\hfil\thepage\hfil}
    \let\@evenfoot\@oddfoot
    \def\@evenhead{\small\@date\hfil\slshape\@pauthor\hfil}
    \def\@oddhead{\small\null\hfil\slshape\@ptitle\hfil}
    \let\@mkboth\@gobbletwo
    \let\sectionmark\@gobble
    \let\subsectionmark\@gobble
   }
\else               
  \def\ps@draft{%
    \def\@oddfoot{\small\@date\hfil\slshape\@pauthor\hfil\upshape\thepage}
    \def\@oddhead{\small\null\hfil\slshape\@ptitle\hfil}
    \let\@mkboth\@gobbletwo
    \let\sectionmark\@gobble
    \let\subsectionmark\@gobble
  }
\fi

\newcommand*{\keywords}[1]{\gdef\@keywords{#1}}
\keywords{}
\newcommand{\keywordsname}{Key words and phrases}

\newcommand*{\subjclass}[1]{\gdef\@subjclass{#1}}
\subjclass{}
\newcommand{\subjclassname}{1991 AMS Mathematics Subject Classification}

\def\@maketitle{%
  \newpage
  \null
  \vskip 2em%
  \begin{center}%
    {\Large\bfseries \@title \par}%
    \vskip 1.5em%
    {\small\scshape
      \lineskip .5em%
      \begin{tabular}[t]{c}%
        \@author
      \end{tabular}\par
    }%
    \vskip 1em%
    {\small\@date}
  \end{center}%
  \par
  \vskip 1.5em
  \begingroup
    \let\@makefnmark\relax \let\@thefnmark\relax
    \ifx\@empty\@subjclass\else
       \@footnotetext{{\itshape\subjclassname}.\enspace\@subjclass.}
    \fi
    \ifx\@empty\@keywords\else
       \@footnotetext{{\itshape\keywordsname}.\enspace\@keywords.}
    \fi
  \endgroup
}

\renewcommand{\part}{\par
   \addvspace{4ex}%
   \@afterindentfalse
   \secdef\@part\@spart}

\def\@part[#1]#2{%
    \ifnum \c@secnumdepth >\m@ne
      \refstepcounter{part}%
      \addcontentsline{toc}{part}{\thepart\hspace{1em}#1}%
    \else
      \addcontentsline{toc}{part}{#1}%
    \fi
    {\parindent \z@ \raggedright
     \interlinepenalty \@M
     \reset@font
     \ifnum \c@secnumdepth >\m@ne
       \large\bfseries \partname~\thepart
       \par\nobreak
     \fi
     \Large \bfseries #2%
     \markboth{}{}\par}%
    \nobreak
    \vskip 3ex
    \@afterheading}
\def\@spart#1{%
    {\parindent \z@ \raggedright
     \interlinepenalty \@M
     \reset@font
     \Large \bfseries #1\par}%
     \nobreak
     \vskip 3ex
     \@afterheading}
\def\@endpart{\vfil\newpage
              \if@twoside
                \hbox{}%
                \thispagestyle{empty}%
                \newpage
              \fi
              \if@tempswa
                \twocolumn
              \fi}
\renewcommand{\section}{\@startsection {section}{1}{\z@}%
                                   {-3.5ex \@plus -1ex \@minus -.2ex}%
                                   {2.3ex \@plus.2ex}%
                                   {\reset@font\large\bfseries}}
\renewcommand{\subsection}{\@startsection{subsection}{2}{\z@}%
                                     {-3.25ex\@plus -1ex \@minus -.2ex}%
                                     {1.5ex \@plus .2ex}%
                                     {\reset@font\normalsize\bfseries}}
\renewcommand{\subsubsection}{\@startsection{subsubsection}{3}{\z@}%
                                     {-3.25ex\@plus -1ex \@minus -.2ex}%
                                     {1.5ex \@plus .2ex}%
                                     {\reset@font\normalsize\bfseries}}
\renewcommand{\paragraph}{\@startsection{paragraph}{4}{\z@}%
                                    {3.25ex \@plus1ex \@minus.2ex}%
                                    {-1em}%
                                    {\reset@font\normalsize\bfseries}}
\renewcommand{\subparagraph}{\@startsection{subparagraph}{5}{\parindent}%
                                       {3.25ex \@plus1ex \@minus .2ex}%
                                       {-1em}%
                                      {\reset@font\normalsize\bfseries}}

\renewcommand{\theenumi}{\alph{enumi}}
\renewcommand{\labelenumi}{(\theenumi)}
\renewcommand{\theenumii}{\roman{enumii}}

\renewcommand{\p@enumii}{\theenumi.}
\renewcommand{\theenumiii}{\Alph{enumiii}}

\renewcommand{\p@enumiii}{\theenumi.\theenumii.}

\renewcommand{\p@enumiv}{\p@enumiii\theenumiii.}

\@addtoreset{equation}{section}

\RequirePackage{amsthm}

\renewenvironment{proof}[1][\proofname]{\par
  \normalfont
  \topsep6\p@\@plus6\p@ \trivlist
  \item[\hskip\labelsep\slshape
    #1\@addpunct{.}]\ignorespaces
}{%
  \qed\endtrivlist
}

\theoremstyle{plain}

\newtheorem{theorem}{Theorem}[section]
\newtheorem{proposition}[theorem]{Proposition}
\newtheorem{lemma}[theorem]{Lemma}
\newtheorem{corollary}[theorem]{Corollary}

\theoremstyle{definition}

\newtheorem{definition}[theorem]{Definition}

\newtheorem{example}[theorem]{Example}

\newtheorem{remark}[theorem]{Remark}

\RequirePackage{amsmath}
\RequirePackage{amssymb}

\newcommand{\C}{\mathbb{C}}
\newcommand{\N}{\mathbb{N}}

\newcommand{\Z}{\mathbb{Z}}

\newcommand{\DSum}{\bigoplus}

\DeclareMathOperator{\im}{im}

\renewcommand{\to}[1][]{\xrightarrow[#1]{}}
\newcommand{\from}[1][]{\xleftarrow[#1]{}}
\newcommand{\isoto}[1][]{\xrightarrow[#1]{\sim}}

\newcommand{\Endo}[1][]{\mathrm{End}_{\raise1.5ex\hbox to.1em{}#1}}

\newcommand{\Hom}[1][]{\mathrm{Hom}_{\raise1.5ex\hbox to.1em{}#1}}

\newcommand{\RHom}[1][]{\mathrm{RHom}_{\raise1.5ex\hbox to.1em{}#1}}

\newcommand{\Ext}[2][]{\mathrm{Ext}_{\raise1.5ex\hbox to.1em{}#1}^{#2}}

\newcommand{\THom}[1][]{\mathrm{THom}_{\raise1.5ex\hbox to.1em{}#1}}

\newcommand{\Tens}[1][]{\mathbin{\otimes_{\raise1.5ex\hbox to-.1em{}#1}}}

\newcommand{\LTens}[1][]{\mathbin{\otimes_{\raise1.5ex\hbox to-.1em{}#1}^{L}}}

\newcommand{\Tor}[2][]{\mathrm{Tor}^{\raise1.5ex\hbox to.1em{}#1}_{#2}}

\def\sha{\mathcal{A}}

\def\she{\mathcal{E}}
\def\shf{\mathcal{F}}
\def\shg{\mathcal{G}}

\def\shi{\mathcal{I}}
\def\shj{\mathcal{J}}

\def\shl{\mathcal{L}}
\def\shm{\mathcal{M}}
\def\shn{\mathcal{N}}

\newcommand{\sect}{\Gamma}

\renewcommand{\hom}[1][]{{\mathcal{H}om}_{\raise1.5ex\hbox to.1em{}#1}}

\newcommand{\rhom}[1][]{{R\mathcal{H}om}_{\raise1.5ex\hbox to.1em{}#1}}

\newcommand{\ext}[2][]{{\mathcal{E}xt}_{\raise1.5ex\hbox to.1em{}#1}^{#2}}

\newcommand{\thom}[1][]{{T\mathcal{H}om}_{\raise1.5ex\hbox to.1em{}#1}}

\newcommand{\tens}[1][]{\mathbin{\otimes_{\raise1.5ex\hbox to-.1em{}#1}}}

\newcommand{\ltens}[1][]{\mathbin{\otimes_{\raise1.5ex\hbox to-.1em{}#1}^{L}}}

\newcommand{\tor}[2][]{{\mathcal{T}or}^{\raise1.5ex\hbox to.1em{}#1}_{#2}}

\newcommand\etens{\mathbin{\boxtimes}}

\newcommand{\opb}[1]{#1^{-1}}

\newcommand{\GHom}[1][]{\mathrm{GHom}_{\raise1.5ex\hbox to.1em{}#1}}

\newcommand{\GExt}[2][]{\mathrm{GExt}_{\raise1.5ex\hbox to.1em{}#1}^{#2}}

\newcommand{\FHom}[1][]{\mathrm{FHom}_{\raise1.5ex\hbox to.1em{}#1}}

\newcommand{\ghom}[1][]{{\mathcal{GH}om}_{\raise1.5ex\hbox to.1em{}#1}}

\newcommand{\gext}[2][]{{\mathcal{GE}xt}_{\raise1.5ex\hbox to.1em{}#1}^{#2}}

\newcommand{\fhom}[1][]{{\mathcal{FH}om}_{\raise1.5ex\hbox to.1em{}#1}}

\newcommand{\gr}{\mathop{\mathcal{G}r}\nolimits}

\newcommand{\f}{\mathcal{F}}

\newcommand{\tenstop}[1][]{\mathbin{\hat{\otimes}_{\raise1.5ex\hbox to-.1em{}#1}}}

\newcommand{\homtop}[1][]{\mathcal{L}_{\raise1.5ex\hbox to.1em{}#1}}

\newcommand{\Homtop}[1][]{\mathrm{L}_{\raise1.5ex\hbox to.1em{}#1}}

\newcommand{\D}{\mathcal{D}}

\newcommand{\E}{\mathcal{E}}

\renewcommand{\O}{\mathcal{O}}

\def\absdoim#1{\underline{#1}_*}
\def\reldoim[#1]#2{\underline{#2}_{|{#1}*}}
\def\doim{\@ifnextchar [{\reldoim}{\absdoim}}

\def\absdeim#1{\underline{#1}_*}
\def\reldeim[#1]#2{\underline{#2}_{|{#1}*}}
\def\deim{\@ifnextchar [{\reldeim}{\absdeim}}

\def\absdopb#1{\underline{#1}^{-1}}
\def\reldopb[#1]#2{\underline{#2}_{|{#1}}^{-1}}
\def\dopb{\@ifnextchar [{\reldopb}{\absdopb}}

\def\absboim#1{\underline{\underline{#1}}_*}
\def\relboim[#1]#2{\underline{\underline{#2}}_{|{#1}*}}
\def\boim{\@ifnextchar [{\relboim}{\absboim}}

\def\absbeim#1{\underline{\underline{#1}}_*}
\def\relbeim[#1]#2{\underline{\underline{#2}}_{|{#1}*}}
\def\beim{\@ifnextchar [{\relbeim}{\absbeim}}

\def\absbopb#1{\underline{\underline{#1}}^*}
\def\relbopb[#1]#2{\underline{\underline{#2}}_{|{#1}}^*}
\def\bopb{\@ifnextchar [{\relbopb}{\absbopb}}

\newcommand{\eu}{\mathop{\mathrm{eu}}\nolimits}

\RequirePackage{xypic}

\newcommand{\eqdot}{\mathbin{:=}}

\newcommand{\genclass}[3]{[#1]_{#3}^{#2}}
\newcommand{\lsclass}[1]{\genclass{#1}{1}{\C^\times}}
\newcommand{\lbclass}[1]{\genclass{#1}{1}{\O^\times}}
\newcommand{\stkclass}[1]{\genclass{#1}{2}{\C^\times}}
\newcommand{\tlbclass}[1]{\genclass{#1}{1}{{\overline\O}^\times}}
\newcommand{\tdoclass}[1]{\genclass{#1}{1}{d\O}}

\newcommand{\OTX}[1][]{\O_{T^*X}\def\temp{#1}\ifx\temp\empty\else(#1)\fi}
\newcommand{\OV}[1][]{\O_{V}\def\temp{#1}\ifx\temp\empty\else(#1)\fi}
\newcommand{\JV}[1][]{\shj_V\def\temp{#1}\ifx\temp\empty\else(#1)\fi}

\newcommand{\filt}[2][]{\mathsf{F}_{#1}#2}
\newcommand{\optionfilt}[2][]{\def\temp{#1}\ifx\temp\empty#2\else\filt[#1]{#2}\fi}
\newcommand{\Vfilt}[3][]{\mathsf{F}^{#2}_{#1}#3}
\newcommand{\optionVfilt}[3][]{\def\temp{#1}\ifx\temp\empty\E_{#2}\else\Vfilt[#1]{#2}{#3}\fi}

\newcommand{\DX}[1][]{\optionfilt[#1]{\D_X}}
\newcommand{\EX}[1][]{\optionfilt[#1]{\E_X}}
\newcommand{\EC}[1][]{\optionfilt[#1]{\E_\C}}
\newcommand{\EV}[1][]{\optionVfilt[#1]{V}{\E_X}}

\newcommand{\IV}[1][]{\shi_V\def\temp{#1}\ifx\temp\empty\else^{(#1)}\fi}
\newcommand{\ITX}[1][]{\shi_{T^*X}\def\temp{#1}\ifx\temp\empty\else^{(#1)}\fi}
\newcommand{\IVdot}[1][]{\shi_{\bigdot V}\def\temp{#1}\ifx\temp\empty\else^{(#1)}\fi}

\renewcommand{\eu}{\operatorname{\mathit{eu}}}
\newcommand{\Lie}[1][]{\operatorname{\mathsf{L}}\def\temp{#1}\ifx\temp\empty\else^{(#1)}\fi}

\newcommand{\catMod}{\mathsf{Mod}}

\newcommand{\shHom}[1][]{{\mathcal{H}om}_{\raise1.5ex\hbox to.1em{}#1}}
\newcommand{\shEndo}[1][]{{\mathcal{E}nd}_{\raise1.5ex\hbox to.1em{}#1}}
\newcommand{\shAut}[1][]{{\mathcal{A}ut}_{\raise1.5ex\hbox to.1em{}#1}}

\newcommand{\stks}{\mathfrak{S}}
\newcommand{\stkt}{\mathfrak{T}}
\newcommand{\stku}{\mathfrak{U}}

\newcommand{\stkHom}[1][]{\mathfrak{Hom}_{\raise1.5ex\hbox to.1em{}#1}}
\newcommand{\stkMod}{\mathfrak{Mod}}

\newcommand{\rk}{\mathop{\mathrm{rk}}}

\newcommand{\twst}{\mathsf{t}}

\newcommand{\DD}{\mathbb{D}}
\renewcommand{\deim}[1]{{\DD #1}_{!}}
\renewcommand{\doim}[1]{{\DD #1}_{*}}
\renewcommand{\dopb}[1]{\DD #1^{*}}
\newcommand{\dzeroopb}[1]{#1^{*}}

\newcommand{\dhom}[1][]{{\DD\mathcal{H}om}_{\raise1.5ex\hbox to.1em{}#1}}

\newcommand{\tdotwist}[2][\D]{{#1}_{#2}}

\newcommand{\CX}{\C_X}

\newcommand{\OX}{\O_X}

\newcommand{\DV}{\D_V}

\renewcommand{\AA}{\mathbb{A}}
\newcommand{\PP}{\mathbb{P}}
\newcommand{\bb}{\PP^*}
\newcommand{\OP}{\O_\PP}
\newcommand{\DP}{\D_\PP}
\newcommand{\GG}{\mathbb{G}}

\newcommand{\TT}{\mathbb{T}}
\newcommand{\CG}{\C_\GG}
\newcommand{\OG}{\O_\GG}
\newcommand{\DG}{\D_\GG}
\newcommand{\FF}{\mathbb{F}}

\newcommand{\stktimes}{\circledast}

\newcommand{\stktens}{\mathbin{\stktimes}}

\newcommand{\stkopb}[1]{#1^{\stktimes}}

\newcommand{\stkinv}[2][-1]{#2^{\stktimes#1}}

\def\@ptitle{Non-existence of twisted wave equations}
\makeatother

\author{Andrea D'Agnolo \and Pierre Schapira}
\title{On twisted microdifferential modules I.\\
Non-existence of twisted wave equations}
\date{}

\begin{document}

\maketitle

\begin{abstract}
Using the notion of subprincipal symbol, we give a necessary condition for the
existence of twisted $\D$-modules simple along a smooth
involutive submanifold of the cotangent bundle to a complex manifold.
As an application, we prove that there are no
generalized massless field equations with non trivial twist on grassmannians, 
and in particular that 
the Penrose transform does not extend to the twisted case.
\end{abstract}

\section*{Introduction}

Let $\TT$ be a $4$-dimensional complex vector space,
$\PP$ the $3$-dimensional projective
space of lines in $\TT$, and $\GG$ the $4$-dimensional 
grassmannian of $2$-planes in $\TT$. 
According to Penrose $\GG$
is a complexification of a conformal compactification of the flat
Minkowski space. Denote by $\shm_{(h)}$ the $\DG$-module associated
with the massless field equations of helicity $h\in\frac12\Z$. 
The Penrose correspondence realizes $\shm_{(1+m/2)}$ as the transform of the 
$\DP$-module associated with the line bundle $\OP(m)$, for $m\in\Z$.
For $\lambda\in\C$, $\OP(\lambda)$ makes sense in the theory of twisted
sheaves. It is then a natural question to ask whether 
the Penrose correspondence extends
to the twisted case. In particular, are there ``massless field equations" 
of complex helicity $h\notin \frac12\Z$?

The $\DG$-modules $\shm_{(h)}$ are simple along a smooth 
involutive submanifold $V$ of the cotangent bundle to $\GG$, 
which is given by the geometry of the integral transform.
In this paper we give a negative answer to the question raised above: 
for topological reasons, there are no simple $\DG$-modules along $V$ 
with non trivial twist. Indeed, this is a corollary of the following more general result. 

Let $X$ be a complex manifold, and $V$ a conic involutive submanifold of its
cotangent bundle. Denote by $\D_{\Omega_{V/X}^{1/2}}$ the 
ring of differential operators on $V$ acting on relative half-forms and by $\D_{\Omega_{V/X}^{1/2}}^{bic}(0)$ the  subring of
operators homogeneous of degree $0$ and commuting with the functions
which are locally constant on the bicharacteristic leaves. 
The ring of microdifferential operators $\EX$ is 
endowed with the so-called $V$-filtration $\{\EV[k]\}_{k\in\Z}$ and 
by a result of Kashiwara-Oshima, there
is a natural isomorphism of rings 
$\EV[0]/\EV[-1] \isoto \D_{\Omega_{V/X}^{1/2}}^{bic}(0)$.

Let $\stks$ be a stack of twisted sheaf on $X$, and consider the category of
twisted microdifferential modules $\catMod(\EX;\stks)$.
One says that a twisted microdifferential module is 
simple along $V$ if it can be endowed with a good $V$-filtration 
whose associated graded module is locally isomorphic to
$\OV[0]$.

Let $\Sigma$ be a smooth bicharacteristic leaf of $V$.
Recall that stacks of twisted sheaves on $X$ are classified by 
$H^2(X;\CX^\times)$, and denote by
$\stkclass\stks$ the class of $\stks$. 
Our main result consists in associating to $\stkclass\stks$ a class 
in $H^2(\Sigma;\C_\Sigma^\times)$ whose triviality
is a necessary condition for the existence of a twisted
microdifferential 
module in $\catMod(\EX;\stks)$ simple along $V$. 

Let us briefly describe our construction. Let $\shm$ be a twisted
microdifferential 
module in $\catMod(\EX;\stks)$ which is simple along $V$.
By definition, $\shm$ has a good $V$-filtration, and we denote by
$\overline\shm$ 
its associated graded module.

\begin{itemize}
\item[(i)]
By Kashiwara-Oshima's result, we consider $\overline\shm$ as an object of $\catMod(\DV^{bic}(0);\stkt)$. Here, $\stkt$ 
is a stack of twisted sheaves on $V$ whose class in 
$H^2(V;\C_V^\times)$ is the product of the pull back of 
$\stkclass\stks$ by the class of the stack containing the inverse relative 
half-forms $\Omega_{V/X}^{-1/2}$. 
\item[(ii)]
The restriction of $\overline\shm$ to $\Sigma$ is a flat connection $\shl$ of rank one in $\catMod(\D_\Sigma; \stku)$, where $\stku$ is a stack of 
twisted sheaves on $\Sigma$ whose class 
$\stkclass\stku\in H^2(\Sigma;\C_\Sigma^\times)$ is 
the restriction of $\stkclass\stkt$. 
\item[(iii)]
By the Riemann-Hilbert correspondence, $\shl$ is associated with a 
local system of rank one on $\stku$.  
Since there are no local systems of rank one with non trivial twist, 
the triviality of $\stkclass\stku$ is a necessary condition
for the existence of a twisted microdifferential 
module in $\catMod(\EX;\stks)$ simple along $V$. 
\end{itemize}

We would like to thank Masaki Kashiwara for extremely useful 
conversations and helpful insights.

\section{Review on twisted sheaves}

In this section we briefly review the notions of twisted sheaves. 
References are made to \cite{K3,KQ}, see also \cite{DP}.

\medskip

Let $X$ be a complex analytic manifold, 
$\OX$ its structure sheaf, and denote by $\CX$ the constant 
sheaf with stalk $\C$. If $\sha$ is a sheaf of $\C$-algebras on $X$, 
we denote by $\catMod(\sha)$ the category of sheaves of $\sha$-modules
on $X$ and by $\stkMod(\sha)$ the corresponding $\C$-stack, 
$U\mapsto \catMod(\sha|_U)$. We denote by $\sha^\times$ the sheaf of
invertible sections of $\sha$.

The short exact sequence of abelian groups
$$
1\to \CX^\times\to \OX^\times \to \OX^\times/\CX^\times \to 1
$$
induces the exact sequence
\begin{equation}
\label{eq:alphabetagamma}
H^1(X;\CX^\times) \to[\alpha] H^1(X;\OX^\times)\to[\beta]
             H^1(X;\OX^\times/\CX^\times)\to[\delta] H^2(X;\CX^\times).
\end{equation}
Note that the isomorphism $d \log\colon\OX^\times/\CX^\times\isoto d\O$
induces an  isomorphism
\begin{equation}
\label{eq:iota2}
\iota \colon H^1(X; \OX^\times/\CX^\times) \isoto H^1(X; d\OX).
\end{equation}
The $\C$-vector space structure of $H^1(X; d\OX)$ thus gives a meaning to $\lambda \cdot c$ for $c\in H^1(X; \OX^\times/\CX^\times)$ and $\lambda\in\C$.

We will consider several characteristic classes with values in these
cohomology groups.

\begin{itemize}
\item[$\bullet$]
A local system is a $\CX$-module locally free of finite rank. 
To a local system $L$ of rank one corresponds a 
class $\lsclass L\in H^1(X;\CX^\times)$ which characterizes $L$ 
up to isomorphisms of $\CX$-modules.

\item[$\bullet$]
A line bundle is an $\OX$-module locally free of rank one.
To a line bundle $\shl$ on $X$ corresponds a
class $\lbclass\shl\in H^1(X;\OX^\times)$ which characterizes $\shl$ 
up to isomorphisms of $\OX$-modules.

\item[$\bullet$]
A stack of twisted sheaves is a $\C$-stack locally $\C$-equivalent to
$\stkMod(\CX)$. To a stack of twisted sheaves $\stks$ corresponds a class 
$\stkclass\stks\in H^2(X;\C^\times_X)$ which
characterizes $\stks$ up to $\C$-equivalences. 
Objects of $\stks(X)$ are called twisted sheaves.
\end{itemize}

Recall that $\stkclass\stks$ has the following description 
using Cech cohomology. 
Let $X=\bigcup_i U_i$ be an open covering 
such that there are $\C$-equivalences 
$\varphi_i\colon \stks|_{U_i}\to\stkMod(\C_{U_i})$. 
By Morita theory, the auto-equivalence
$\varphi_i\circ\varphi_j^{-1}$ of $\stkMod(\C_{U_{ij}})$ are 
given by $G\mapsto G\tens L_{ij}$ for a local system $L_{ij}$ of rank
one. By refining the covering we may assume that $L_{ij} \simeq\C_{U_{ij}}$. 
The isomorphisms $L_{ij} \tens L_{jk} \simeq L_{ik}$ on $U_{ijk}$ are 
then multiplication by locally constant functions  $c_{ijk}\in\sect(U_{ijk}; \C_X^\times)$. 
The class $\stkclass\stks$ is described by the Cech cocycle $\{c_{ijk}\}$.

A twisted sheaf $F\in\stks(X)$ is described by a family of 
sheaves $F_i\in\catMod(\C_{U_i})$ and isomorphisms 
$\theta_{ij} \colon F_j|_{U_{ij}} \to F_i|_{U_{ij}}$ 
satisfying $\theta_{ij} \circ \theta_{jk} = c_{ijk} \theta_{ik}$.

Let $\stks$ be a stack of twisted sheaves on $X$ and let $\sha$ be a sheaf 
of $\C$-algebras on $X$.
We denote by $\stkMod(\sha;\stks)$ the stack of left 
$\sha$-modules in $\stks$. 

\begin{itemize}
\item[$\bullet$]
A twisted line bundle is a pair $(\stks,\shf)$ of a stack of twisted sheaves 
$\stks$ and an object $\shf\in \catMod(\OX;\stks)$ locally 
free of rank one over $\OX$. To a twisted line bundle corresponds a class 
$\tlbclass{\stks,\shf}\in H^1(X;\OX^\times/\CX^\times)$ which
characterizes it up to the following equivalence relation: two twisted
line bundles $(\stks,\shf)$ and $(\stkt,\shg)$ are equivalent if there exist a
$\C$-equivalence $\varphi\colon \stks\to\stkt$ and an
isomorphism $\varphi(\shf) \simeq \shg$ in $\stkt(X)$. 
\end{itemize}

Let $(\stks,\shf)$ be a twisted line bundle and let $X=\bigcup_i U_i$
be an open covering such that there are $\C$-equivalences
$\varphi_i\colon \stks|_{U_i}\to\stkMod(\C_{U_i})$, and denote by
$\{c_{ijk}\}$ the Cech cocycle of $\stkclass\stks$. 
These induce equivalences 
$\varphi_i\colon\stkMod(\O_{U_i};\stks|_{U_i})\to\stkMod(\O_{U_i})$
and $\shf$ is described by a family of line bundles
$\shf_i\in\catMod(\O_{U_i})$ and isomorphisms $\theta_{ij} \colon
\shf_j|_{U_{ij}} \to \shf_i|_{U_{ij}}$. By refining the covering,
we may assume that there are nowhere vanishing sections
$s_i\in\sect(U_i; \shf_i)$, so that $\shf_i\simeq \O_{U_j}$. Hence
$\theta_{ij}$ are multiplications by the sections $f_{ij} =s_i/\theta_{ij}(s_j)\in\sect(U_{ij};\OX^\times)$, so that $f_{ij} f_{jk} = c_{ijk}f_{ik}$. 
The class $\tlbclass{\stks,\shf}$ is thus
described by the Cech hyper-cocycle $\{f_{ij}, c_{ijk}\}$.

The characteristic classes constructed above are related (up to sign) as follows,
using the exact sequence \eqref{eq:alphabetagamma}:
\begin{enumerate}
\item
if $L$ is a local system of rank one, then
$\alpha(\lsclass L)=\lbclass{L\tens\OX}$,
\item
if $\shl$ is a line bundle, then
$\beta(\lbclass\shl)=\tlbclass{\stkMod(\CX),\shl}$,
\item
if $(\stks,\shf)$ is a twisted line bundle,
then $\delta(\tlbclass{\stks,\shf})=\stkclass\stks$. 
\end{enumerate}

The next result will play an essential role in the proof of Theorem \ref{th:main}.
It immediately follows from the Morita theory for stacks.

\begin{proposition}\label{pr:trivialstk}
A stack of twisted sheaves $\stks$ is globally $\C$-equivalent to $\stkMod(\CX)$
if and only if there exists an object $F\in\stks(X)$ locally free of rank one over $\C$. 
\end{proposition}

\begin{example}
For $\shl$ an untwisted line bundle, and $\lambda\in\C$,
there is a twisted line bundle $(\stks_{\shl^\lambda},\shl^\lambda)$
whose class $\tlbclass{\stks_{\shl^\lambda},\shl^\lambda}$ is described as follows.
Let $X=\bigcup_i U_i$ be an open covering such that there are nowhere 
vanishing sections $s_i\in\sect(U_i;\shl)$, and set $g_{ij} = s_i/s_j$. 
Choose a determination $f_{ij}$ for the ramified function $g_{ij}^\lambda$ 
on $U_{ij}$. Then $f_{ij} f_{jk}$ and $f_{ik}$ are different
determinations of $g_{ik}^\lambda$, so that $f_{ij} f_{jk} =c_{ijk}f_{ik}$ 
for some $c_{ijk}\in \sect(U_{ijk};\CX^\times)$. Then $\tlbclass{\stks_{\shl^\lambda},\shl^\lambda}$ is described by the Cech hyper-cocycle $\{f_{ij}, c_{ijk}\}$. Since $d\log f_{ij} = \lambda d\log g_{ij}$, we have
$$
\tlbclass{\stks_{\shl^\lambda},\shl^\lambda} = \lambda \cdot \beta(\lbclass\shl)\quad \text{in } H^1(X;\OX^\times/\CX^\times), 
$$
where the action of $\lambda$ on $\beta(\lbclass\shl)$ is induced by the isomorphism \eqref{eq:iota2}.

Note that $\shl^\lambda$ is unique up to tensoring by a local system of rank one. 
\end{example}

\subsection*{Operations on stacks}
Consider two stacks $\stks$ and $\stks'$
of twisted sheaves on $X$ (here, $X$ is simply a topological space, or
even a site).
There are stacks of twisted sheaves $\stks \stktens \stks'$ and 
$\stkinv\stks$ on $X$ such that if $F\in\stks(X)$ and $F'\in\stks'(X)$ 
are twisted sheaves, then 
$F\tens F' \in (\stks \stktens \stks')(X)$ and if $F$ is a local
system of rank one, then $F^{-1}\in \stkinv\stks$. Moreover, 
\begin{align*}
\stkclass{\stks \stktens \stks'} &= \stkclass{\stks} \cdot \stkclass{\stks'}\\
\stkclass{\stkinv\stks} &= (\stkclass{\stks})^{-1}.
\end{align*}
If $f\colon Y\to X$ is a morphism of topological spaces (or of sites), 
there exists a stack of twisted sheaves $\stkopb f\stks$ on $Y$ such
that if $F\in\stks(X)$, then
$\opb f F \in (\stkopb f\stks)(Y)$.
Moreover,
$$
\stkclass{\stkopb f\stks} = f^\sharp(\stkclass{\stks}).
$$
Here, for $\twst, \twst'\in H^2(X;\C^\times_X)$, we denote by 
$\twst \cdot \twst'$ and $\twst^{-1}$ the product and the inverse in $H^2(X;\C^\times_X)$, 
respectively, and by
$f^\sharp\twst\in H^2(Y;\C^\times_Y)$ the pull-back.

Let $(\stks_\shf,\shf)$ and $(\stks_\shg,\shg)$ be twisted
line bundles on $X$, and consider the associated twisted line bundles 
$(\stks_{\shf^{-1}},\shf^{-1})$ and $(\stks_{\shf\tens[\O]\shg},\shf\tens[\O]\shg)$ on $X$, and $(\stks_{\f^*\shf},f^*\shf)$ on $Y$. Then there are $\C$-equivalences
\begin{align*}
\stks_{\shf^{-1}} &\simeq \stkinv{\stks_\shf}, \\
\stks_{\shf\tens[\O]\shg} &\simeq \stks_\shf \stktens \stks_\shg, \\
\stks_{f^*\shf} &\simeq \stkopb f \stks_\shf.
\end{align*}

\section{Review on twisted differential operators}

In this section we briefly review the notions of
twisted differential operators. References are made
to \cite{K2,BB} (see also \cite{DP} for an exposition).

Recall that $X$ denotes a complex analytic manifold 
and $\DX$ the sheaf of finite order differential operators on $X$. 
Recall that an automorphism of $\DX$ as an $\OX$-ring is described by
a closed one-form.

\begin{itemize}

\item[$\bullet$]
A ring of twisted differential operators (a t.d.o.\ ring  for short) 
is a sheaf of $\OX$-rings locally isomorphic to $\DX$. To a t.d.o.\
ring $\sha$ 
corresponds a class $\tdoclass\sha\in H^1(X;d\OX)$ which characterizes 
$\sha$ up to isomorphism of $\OX$-rings.

\end{itemize}

Let $(\stks,\shf)$ be a twisted line bundle. 
An example of t.d.o.\ ring is given by
$$
\D_\shf = \shf\tens[\O]\DX\tens[\O]\shf^{-1},
$$
where $\shf^{-1} = \hom[\O](\shf,\OX)$. 
Notice that $\shf^{-1} \in \catMod(\OX;\stkinv\stks)$, 
so that $\D_\shf$ is untwisted as a sheaf.

Let $\{f_{ij}, c_{ijk}\}$ be a Cech hyper-cocycle for $\tlbclass{\stks,\shf}$ 
attached to the covering  $X=\bigcup_i U_i$, where $f_{ij} = s_i/\theta_{ij}(s_j)$ 
for $s_i\in\sect(U_i;\shf_i)$. The sections of $\D_\shf$ are described by 
families $s_i \otimes P_i \otimes s_i^{-1}$, where $P_i\in\sect(U_i;\DX)$ and
\begin{equation}
\label{eq:shashf}
P_i = f_{ji} \cdot P_j \cdot f_{ij} \quad \text{in } \sect(U_{ij};\DX).
\end{equation}
The isomorphism $\iota$ in \eqref{eq:iota2} is then described  by $\iota(\tlbclass{\stks,\shf}) = \tdoclass{\D_\shf}$.
In particular, to any t.d.o.\ ring $\sha$ is associated a twisted line bundle $\shf$, unique up to tensoring by a local system of rank one, 
such that 
$\sha \simeq \D_\shf$ as an $\OX$-ring. 

Let $(\stks,\shf)$ be a twisted line 
bundle and $\stkt$ a stack of twisted sheaves on $X$.
There is an equivalence of $\C$-stacks
\begin{align}\label{eq:DF-SeqvD-FS}
\stkMod(\D_{\shf};\stkt) &\to\stkMod(\D_X;\stkinv\stks\stktens\stkt) \\
\shm &\mapsto \shf^{-1} \tens[\O] \shm.\nonumber
\end{align}

Denote by $\Theta_X$ the sheaf of vector fields and by $\Omega_X$ the sheaf 
of forms of maximal degree. 
We end this section by giving an explicit description, which will be of use later on, of the t.d.o.\ ring 
$\tdotwist{\Omega_X^\lambda}$ for $\lambda\in\C$. Let $v\in\Theta_X$. 
Recall that the Lie derivative $\Lie(v)$ acts on differential forms of
any degree, in particular on $\OX$, where $\Lie(v)(a) = v(a)$, 
and on $\Omega_X$. 
Let $\omega$ be a nowhere vanishing local section of $\Omega_X$. 
One checks that the morphism
\begin{align}
\label{eq:Lielambda}
\Lie[\lambda] \colon \Theta_X &\to \tdotwist{\Omega_X^\lambda}  = 
\Omega_X^\lambda \tens[\O] \DX \tens[\O] \Omega_X^{-\lambda} \\
\notag 
v & \mapsto \omega^\lambda \tens (v + \lambda \frac{\Lie(v)(\omega)}{\omega}) 
\tens \omega^{-\lambda}
\end{align}
is well defined and independent from the choice of $\omega$. 
(Here $\Lie(v)(\omega)/\omega = a$, 
where $a\in \OX$ is such that $\Lie(v)(\omega) = a\omega$.)
Then $\tdotwist{\Omega_X^\lambda}$ is generated by $\OX$ 
and $\Lie[\lambda](\Theta_X)$ 
with the relations
\begin{eqnarray}
\Lie[\lambda](av) &=& a \cdot \Lie[\lambda](v) + \lambda v(a), \\
{[} \Lie[\lambda](v) , a {]} &=& v(a), \\
{[} \Lie[\lambda](v) , \Lie[\lambda](w) {]} &=& \Lie[\lambda]([v,w]), 
\end{eqnarray}
for $a\in\OX$, and $v,w\in\Theta_X$. Of course, $\Lie[0](v) = v$ 
and $\Lie[1](v) = \Lie(v)$.

\section{Microdifferential operators on involutive submanifolds}
\label{se:microdiff}

In this section we recall the notions of microdifferential operators 
and $V$-filtration. 
References are made to \cite{SKK,KO} (see also~\cite{K1,K3,S} for 
expositions).

\medskip
Let $W$ be a complex manifold. In this paper, 
by a submanifold of $W$, we mean a 
smooth locally closed complex submanifold. 

Let $X$ be a complex manifold, and denote by $\pi\colon T^*X\to X$ its 
cotangent bundle. Identifying $X$ with the zero-section of $T^*X$, 
one sets $\dot T^*X = T^*X\setminus X$. 

The canonical $1$-form $\alpha_X$ induces a homogeneous symplectic 
structure on $T^*X$. Denote by $\{f,g\}\in\OTX$ the Poisson 
bracket of two functions 
$f,g\in\OTX$ and by
$$
H \colon T^*T^*X \isoto TT^*X
$$
the Hamiltonian isomorphism. 
For $k\in\Z$, denote by $\OTX[k]  \subset \OTX$ the 
subsheaf of functions $\varphi$ homogeneous 
of order $k$, that is, 
satisfying  $\eu(\varphi)=k\cdot\varphi$. 
Here, $\eu = - H(\alpha_X)$ denotes the Euler vector field on $T^*X$, 
the infinitesimal generator of the action of $\C^\times$.

Denote by $\EX$ the ring of microdifferential operators on $T^*X$. 
It is endowed with the order filtration $\{ \EX[m] \}_{m\in\Z}$, 
where $\EX[m]$ is the subsheaf of microdifferential operators of order at most 
$m$. There is a canonical morphism
$$
\sigma_m \colon \EX[m] \to \OTX[m]
$$
called the principal symbol of order $m$.
This morphism induces an isomorphism of 
graded rings 
$\gr\EX = \simeq \DSum\nolimits_k \OTX[k]$.
If $P\in\EX[m]$, $Q\in\EX[l]$, one has
\begin{eqnarray}
\label{eq:sigmaPQ}
\sigma_{m+l}(PQ) &=& \sigma_{m}(P) \sigma_{l}(Q), \\
\label{eq:sigmaPQQP}
\sigma_{m+l-1}([P,Q]) &=& \{ \sigma_{m}(P), \sigma_{l}(Q) \}.
\end{eqnarray}

Let $V \subset T^*X$ be a submanifold and denote by $\JV\subset\OTX$ 
its annihilating ideal.
Recall that $V$ is called homogeneous, or conic, 
 if $\eu\JV \subset \JV$. In this case, 
$\eu_V \eqdot \eu|_V$ is tangent to $V$, and one defines 
$\OV[k] \subset\OV$ similarly to $\OTX[k] \subset \OTX$.
A conic submanifold $V \subset T^*X$ is called involutive if 
for any pair $f,g\in\JV$ of holomorphic functions 
vanishing on $V$, the Poisson bracket $\{f,g\}$ vanishes on $V$.  
A conic involutive submanifold $V$ 
is called regular if $\alpha_X|_V$ never vanishes.

Let $V \subset T^*X$ be a conic involutive submanifold, and set
$$
\IV = \{P\in\EX[1]|_V; \sigma_{1}(P)|_V = 0 \} \subset \EX|_V.
$$
Note that $[\IV,\IV] \subset \IV$.

\begin{definition} \label{de:EV}
Let $V \subset \dot T^*X$ be a conic involutive submanifold. 
One denotes by $\EV$ the subring of 
$\EX|_V$ generated by $\IV$, and one sets $\EV[m] \eqdot \EX[m]|_V \cdot \EV$.
\end{definition}

One easily checks that
$\EV[m] =  \EV \cdot \EX[m]|_V$, and $\EV[m] \cdot \EV[l] \subset \EV[m+l]$.
In particular, $\{ \EV[k] \}_{k\in\Z}$ is an exhaustive filtration of 
$\EX|_V$, called the 
$V$-filtration, and $\EV[-1]$ is a two-sided ideal of $\EV = \EV[0]$.

\begin{example}\label{ex:canV}
Let $(x) = (x_1,\dots,x_n)$ be a local coordinate system on $X$ and
denote by $(x;\xi) = (x_1,\dots,x_n;\xi_1,\dots,\xi_n)$ the associated
homogeneous symplectic local coordinate system on $T^*X$. Recall that
locally, any conic regular involutive submanifold $V$ of codimension $d$ may be
written after a homogeneous symplectic transformation as:
$$
V = \{(x;\xi); \xi_1 = \cdots = \xi_d = 0\}.
$$
In such a case, 
$$
\EV[m] \simeq (\EX[m]|_V)[\partial_{x_1},\dots,\partial_{x_d}].
$$
\end{example}

\section{Systems with simple characteristics}\label{section:simple1}

In this section we recall the notion of systems with simple characteristics.
References are made to \cite{SKK,KO}.  See also \cite{S,K3} for an 
exposition.

\begin{definition}
Let $\shm$ be a coherent $\EX$-module. 
A lattice in $\shm$ is a coherent $\EX[0]$-submodule $\shm_0$ which generates
$\shm$ over $\EX$.
\end{definition}

Recall that if an $\EX[0]$-submodule $\shm_0$ of $\shm$ defined on an
open subset of $\dot T^*X$ is locally of
finite type, then it is coherent.
A lattice $\shm_0$ endows $\shm$ with the filtration
$$
\filt[k]\shm =  \EX[k] \cdot \shm_0.
$$
If $\shm$ is endowed with a filtration $\{\filt[k]\shm\}_k$, its 
associated symbol module is given by
$$
\widetilde{\gr}(\shm) \eqdot \O_{T^*X}\tens[\gr(\EX)]\gr(\shm),
$$
where $\gr(\shm)=\oplus_{k\in\Z}(\filt[k]\shm/\filt[k-1]\shm)$.

\begin{definition}
\label{def:simplemod}
Let $V \subset \dot T^*X$ be a conic involutive submanifold.
\begin{enumerate}
\renewcommand{\theenumi}{\alph{enumi}}
\renewcommand{\labelenumi}{(\theenumi)}
\item
A coherent $\EX$-module 
$\shm$ is simple along $V$ if it is locally generated by
a section $u\in\shm$, called a simple generator, such that denoting by
$\shi_u$ the annihilator ideal of $u$ in $\EX$, the symbol ideal
$\widetilde\gr(\shi_u)$ is
reduced and coincides with the annihilator ideal $\JV$ of $V$ in $\O_{T^*X}$.
\item
A coherent $\EX$-module 
$\shm$ is globally simple along $V$ if it admits a lattice
$\shm_0$  such that 
$\E_V\shm_0\subset\shm_0$ and
$\shm_0/\filt[-1]\shm$ 
is locally isomorphic to $\OV(0)$.
Such an $\shm_0$ is called a $V$-lattice in $\shm$.
\end{enumerate}
\end{definition}

\begin{lemma}\label{le:simplemod}
If $\shm$ is globally simple, then it is simple. 
\end{lemma}

\begin{proof}
Let $\shm_0$ be a $V$-lattice. Choose a local section $u\in\shm_0$ whose image
in  $\shm_0/\filt[-1]\shm$ is a generator
of $\OV(0)$. Then $\shm_0=\EX[0] u+\filt[-1]\shm$ and it follows
that for all $k\leq 0$ 
$$
\shm_0=\EX[0] u+\filt[k]\shm
$$
Since the filtration on $\shm$ is separated (see~\cite{SKK}), $u$
generates $\shm_0$ over $\EX[0]$.
\end{proof}

Let $(t)\in\C$ be a coordinate, and denote by $(t;\tau)\in T^*\C$ the
associated homogeneous symplectic coordinate system. 
Let $V\subset T^*X$ be a conic involutive submanifold, non necessarily regular.
The trick of the dummy variable consists in replacing $V$ with the
conic involutive submanifold $\widetilde V = V\times \dot T^*\C$, 
which is regular.
Let $p\in V$ and $q\in \dot T^*\C$.
If $\Sigma$ is the bicharacteristic leaf of $V$ through $p$, then 
$\Sigma \times \{q\}$
is the bicharacteristic leaf of $\widetilde V$ through $(p,q)$. 

\begin{proposition}
\label{pr:dummy}
Let $\shm$ be a globally simple $\EX$-module along $V$. Then 
$\widetilde \shm = \E_{X\times\C}\tens[\EX\etens\E_\C](\shm\etens\E_\C)$ 
is globally simple along $\widetilde V$.
\end{proposition}

\begin{proof}
Let $\shm_0$ be a $V$-lattice in $\shm$, and set
$$
\widetilde \shm_0 = \filt[0]{\E_{X\times \C}} \tens[{\EX[0]\etens\filt[0]{\E_\C}}](\shm_0 \etens \EC[0]).
$$
Clearly, $\widetilde \shm_0$ is a lattice in $\widetilde \shm$
and moreover, 
$\she_{\widetilde V}\widetilde \shm_0\subset\widetilde\shm_0$.
Note that 
$$
\filt[-1]{\widetilde \shm} = \filt[0]{\E_{X\times \C}} \tens[{\EX[0]\etens\filt[0]{\E_\C}}](\filt[-1]\shm\etens\filt[0]{\E_\C}+\shm_0\etens\filt[-1]{\E_\C}).
$$
Consider the commutative exact diagram of $\EX[0]\etens\filt[0]{\E_\C}$-modules:
{\small
$$
\xymatrix@C=.7em @R=1em{
       &0\ar[d]&0\ar[d]&0\ar[d]&\\
0\ar[r]&{\filt[-1]\shm\etens\filt[-1]{\E_\C}}\ar[r]\ar[d]
               &{\shm_{0}\etens\filt[-1]{\E_\C}}\ar[r]\ar[d]
                       &{(\shm_{0}/\filt[-1]\shm)\etens\filt[-1]{\E_\C}}\ar[r]\ar[d]
                               &0\\
0\ar[r]&{\filt[-1]\shm\etens\filt[0]{\E_\C}}\ar[r]\ar[d]
               &{\shm_{0}\etens\filt[0]{\E_\C}}\ar[r]\ar[d]
                       &{(\shm_{0}/\filt[-1]\shm)\etens\filt[0]{\E_\C}}\ar[r]\ar[d]
                               &0\\
0\ar[r]&{\filt[-1]\shm\etens(\filt[0]{\E_\C}/\filt[-1]{\E_\C})} \ar[r]\ar[d]
               &{\shm_{0}\etens(\filt[0]{\E_\C}/\filt[-1]{\E_\C})} \ar[r]\ar[d]
&{(\shm_{0}/\filt[-1]\shm)\etens(\filt[0]{\E_\C}/\filt[-1]{\E_\C})} \ar[r]\ar[d]
                               &0\\
       &0      &0      &0      &
}
$$
}
It follows that that the sequence 
$$
0\to \filt[-1]\shm\etens\filt[0]{\E_\C}+\shm_0\etens\filt[-1]{\E_\C}
\to\shm_0\etens\filt[0]{\E_\C}\to\shm_0/\filt[-1]\shm\etens\filt[0]{\E_\C}/\filt[-1]{\E_\C}\to 0
$$
is exact.
Since $\filt[0]{\E_{X\times \C}}$ is flat over $\EX[0]\etens\filt[0]{\E_\C}$, we locally have
\begin{align*}
\filt[0]{\widetilde \shm}/\filt[-1]{\widetilde \shm}
&\simeq \filt[0]{\E_{X\times \C}} \tens[{\EX[0]\etens\filt[0]{\E_\C}}](\shm_0/\filt[-1]\shm\etens\filt[0]{\E_\C}/\filt[-1]{\E_\C}) \\
&\simeq\filt[0]{\E_{X\times \C}} \tens[{\EX[0]\etens\filt[0]{\E_\C}}](\OV[0]\etens\O_{\dot T^*\C}(0)) \\
&\simeq \O_{\widetilde V}(0).
\end{align*}
\end{proof}
\begin{remark}
Let $\stks$ be a $\C$-stack of twisted  sheaves on $X$. Then
Definition \ref{def:simplemod}, Lemma \ref{le:simplemod} and 
Proposition \ref{pr:dummy} extend to objects of  
$\catMod(\EX;\stkopb{\pi_X}\stks)$.
\end{remark}

\section{Differential operators on involutive submanifolds}

We recall here the construction of the ring of homogeneous twisted 
differential operators invariant by the bicharacteristic flow.

\medskip

Let $V \subset T^*X$ be a conic regular involutive submanifold and 
denote by $TV^\bot \subset TV$ the symplectic orthogonal to $TV$.
Denote by $\Theta_{V}^\bot \subset \Theta_V$ the sheaf of sections 
of the bundle 
$TV^\bot \to V$, and let
\begin{align*}
\OV^{bic} &\eqdot \{ a\in\OV; v(a)=0 \text{ for any } v\in \Theta_{V}^\bot \},\\
\OV^{bic}(k) &\eqdot \OV^{bic}\cap\OV(k).
\end{align*}
Then $\OV^{bic}$ is the sheaf of holomorphic functions 
locally constant along the 
bicharacteristic leaves of $V$. Consider the ring
$$
\DV^{bic} = \{ P\in\DV; [a,P]=0 \text{ for any } a\in \OV^{bic} \},
$$
and the subring of operators homogeneous of degree zero
$$
\DV^{bic}(0) = \{ P\in\DV^{bic}; [\eu_V, P]=0 \}.
$$ 

\begin{example}\label{ex:canV2}
Let  $(x;\xi) = (x_1,\dots,x_n;\xi_1,\dots,\xi_n)$ be a local 
homogeneous symplectic coordinate system on $T^*X$ and assume that
$$
V = \{(x;\xi); \xi_1 = \cdots = \xi_d = 0\}.
$$
Set $x' = (x_1,\dots,x_d)$, $x'' = (x_{d+1},\dots x_n)$, 
and similarly set $\xi = (\xi',\xi'')$.
One has $(x',x'',\xi'')\in V$, and the bicharacteristic 
leaves of $V$ are the submanifolds  defined by 
$$
\Sigma=\{(x',x'';\xi'');(x'';\xi'')=(x''_0;\xi''_0)\}.
$$ 
The Euler field $\eu_V$ is given by
$$
\eu_V=\sum_{d+1}^n \xi_i\partial_{\xi_i}=\xi''\partial_{\xi''}.
$$
Hence a function locally  constant along the bicharacteristic
leaves depends only on $(x'',\xi'')$. 
A section of $\OV[0]$ is a holomorphic functions in the variable 
$(x',x'',\xi'')$, homogeneous of degree $0$ with respect to $\xi''$.
Moreover a section of
$\D^{bic}_V(0)$ is uniquely written as a finite sum 
\begin{equation}
\label{eq:sectDV}
\sum_{\alpha\in\N^d}a_\alpha \partial_{x'}^\alpha,\mbox{ with }
a_\alpha\in\OV[0]. 
\end{equation} 
\end{example}

Assume that $V$ is regular, and let $j_\Sigma\colon \Sigma\to V$ be the embedding of a bicharacteristic leaf.
Denote by $\shj_\Sigma^{bic}(0)$ the annihilator ideal of 
$\Sigma$ in $\OV^{bic}(0)$,
and note that $\O_\Sigma \simeq \OV^{bic}(0)/\shj_\Sigma^{bic}(0)|_\Sigma$.
Since $\OV^{bic}(0)$ is in the center of $\DV^{bic}(0)$, 
there is a restriction map
\begin{align*}
\dzeroopb{j_\Sigma}\colon 
\catMod(\DV^{bic}(0)) &\to  
\catMod(\D_\Sigma) \\
\shm &\mapsto \O_\Sigma\tens[\OV^{bic}(0)|_\Sigma]\shm|_\Sigma.
\end{align*}
We will be 
interested in the twisted analogue of the above construction. Namely, set
\begin{eqnarray*}
\D_{\Omega_V^{1/2}}^{bic} &\eqdot&
\{ P\in\D_{\Omega_V^{1/2}}; [a,P]=0 \text{ for any } a\in \OV^{bic} \}, \\
\D_{\Omega_V^{1/2}}^{bic}(0) &\eqdot& 
\{ P\in\D_{\Omega_V^{1/2}}^{bic}; [\Lie^{(1/2)}(\eu_V), P]=0 \}.
\end{eqnarray*}
For $p\in\Sigma$, the quotient $T_pV/T_p\Sigma\simeq T_pV/T_pV^\bot$ 
is a symplectic space.
Hence $j_\Sigma^* \Omega_V \simeq \Omega_\Sigma$.
Thus, there is a restriction morphism
\begin{equation}
\label{eq:j*Sigma}
\dzeroopb{j_\Sigma}\colon 
\catMod(\D^{bic}_{\Omega_V^{1/2}}(0)) \to  
\catMod(\D_{\Omega_\Sigma^{1/2}}).
\end{equation}

\section{Subprincipal symbol}
In this section we recall the notion of subprincipal symbol,
and prove the regular involutive analogue
of an isomorphism obtained in~\cite[Lemma 1.5.1]{KK} 
for the Lagrangian case.
References are made to \cite{KO,KK,K1,K3}.

\medskip

As we will recall, the subprincipal symbol is intrinsically defined for
microdifferential operators twisted by half-forms. We will thus consider here the ring
$$
\E_{\Omega_X^{1/2}} = \opb{\pi}\Omega_X^{1/2}\tens[\opb{\pi}\O]\E_X
\tens[\opb{\pi}\O]\opb{\pi}\Omega_X^{-1/2},
$$
instead of $\EX$. All the notions recalled in
Section~\ref{se:microdiff} 
extend to this ring. In particular, its $V$-filtration is defined by
\begin{equation}
\label{eq:twistE}
\begin{cases}
\IV^{\Omega_X^{1/2}}=\{P\in\filt[1]{\E_{\Omega_X^{1/2}}}|_V; \sigma_{1}(P)|_V= 0\}
\simeq 
\opb{\pi}{\Omega_X^{1/2}}\tens[\opb{\pi}\O]\IV\tens[\opb{\pi}\O]\opb{\pi}{\Omega_X^{-1/2}}, \\
\Vfilt[m]V{\E_{\Omega_X^{1/2}}} = 
\opb{\pi}{\Omega_X^{1/2}}\tens[\opb{\pi}\O]\EV[m]\tens[\opb{\pi}\O]
                                                    \opb{\pi}{\Omega_X^{-1/2}},\\
\E_{V,{\Omega_X^{1/2}}} = \Vfilt[0]V{\E_{\Omega_X^{1/2}}}.
\end{cases}
\end{equation}

Let $(x)$ be a local coordinate system on $X$, and denote by 
$(x;\xi)$ the associated  homogeneous symplectic coordinate  system 
on $T^*X$. A
microdifferential operator $P \in \filt[m]{\E_{\Omega_X^{1/2}}}$ 
is then described by its total 
symbol $\{p_k(x;\xi) \}_{k\leq m}$, where $p_k\in\OTX[k]$.
The functions $p_k$ depend on the local 
coordinate system $(x)$ on $X$, except the top degree term $p_m =\sigma_m(P)$
which does not. Recall that the subprincipal symbol
$$
\sigma'_{m-1} \colon \filt[m]{\E_{\Omega_X^{1/2}}} \to \OTX[m-1]
$$
given by 
$$
\sigma'_{m-1}((dx)^{1/2} \tens P \tens (dx)^{-1/2}) =
p_{m-1}(x,\xi) - \frac 12 \sum_i \partial_{x_i}\partial_{\xi_i} p_{m}(x,\xi),
$$
does not depend on the local coordinate system $(x)$ on $X$.
For $P\in\filt[m]{ \E_{\Omega_X^{1/2}} }$, $Q\in\filt[l]{ \E_{\Omega_X^{1/2}} }$, one has
\begin{eqnarray}
\label{eq:subprincipalPQ}
\sigma'_{m+l-1}(PQ) &=& \sigma_{m}(P) \sigma'_{l-1}(Q) +
\sigma'_{m-1}(P) \sigma_{l}(Q) +
\frac12 \{ \sigma_{m}(P), \sigma_{l}(Q) \} , \\
\label{eq:subprincipalPQQP}
\sigma'_{m+l-2}([P,Q]) &=& \{ \sigma_{m}(P), \sigma'_{l-1}(Q) \}
+ \{ \sigma'_{m-1}(P), \sigma_{l}(Q) \} . 
\end{eqnarray}

Let $V\subset T^*X$ be a conic involutive submanifold. 
For $f\in\OTX$, denote by $H_f = H(df) \in TT^*X$ 
its Hamiltonian vector field. Recall that $H$ induces an isomorphism
\begin{equation}
\label{eq:HTVbot}
H \colon T^*_VT^*X \isoto TV^\bot.
\end{equation}
In particular, $H_f|_V$ is tangent to $V$ for $f\in\JV$.
With notations \eqref{eq:Lielambda}, set
\begin{align}
\label{eq:LV}
\shl^0_V \colon \IV^{\Omega_X^{1/2}} &\to \filt[1]{\tdotwist{\Omega_{V}^{1/2}}}, \\
\notag
P &\mapsto \Lie[1/2](H_{\sigma_{1}(P)}|_V) +\sigma'_0(P)|_V .
\end{align}
Using the above relations, one checks that the morphism $\shl^0_V$ does not depend on the choice of coordinates, and satisfies the relation
\begin{align*}
\shl^0_V(AP) &= \sigma_0(A) \shl^0_V(P) , \\
\shl^0_V(PA) &= \shl^0_V(P) \sigma_0(A) , \\
\shl^0_V([P,Q]) &= [ \shl^0_V(P) , \shl^0_V(Q) ],
\end{align*}
for $P,Q\in\IV^{\Omega_X^{1/2}}$ and $A \in \filt[0]{\E_{\Omega_X^{1/2}}}$ (see \cite[\S2]{KO} or \cite[\S8.3]{K3}).
It follows that $\shl^0_V$ extends as a ring morphism
\begin{equation}
\label{eq:LVtemp}
\shl_V\colon \E_{V,\Omega_X^{1/2}} \to \tdotwist{\Omega_{V}^{1/2}}
\end{equation}
by setting $\shl_V(P_1\cdots P_r) = \shl_V^0(P_1) \cdots \shl_V^0(P_r)$, 
for $P_i \in \IV^{\Omega_X^{1/2}}$. 

\begin{theorem}
\label{th:LV1}
Let $V\subset \dot T^*X$ be a conic regular involutive submanifold.
The morphism \eqref{eq:LVtemp} induces a ring isomorphism 
\begin{equation}
\label{eq:LV3}
\shl_V\colon \E_{V,\Omega_X^{1/2}}/\Vfilt[-1]V{\E_{\Omega_X^{1/2}}} \isoto 
\D^{bic}_{\Omega_V^{1/2}}(0).
\end{equation}
\end{theorem}

It is possible to show that the above statement holds even 
without the assumption of regularity for $V$ (for example, 
the Lagrangian case is obtained in~\cite[Lemma~1.5.1]{KK}).

\begin{proof}
The statement is local. We may thus assume that 
$\Omega_X \simeq \OX$ and $\Omega_V \simeq \OV$, 
so that we are reduced to prove the isomorphism
$$
\shl_V\colon \EV/\EV[-1] \isoto 
\DV^{bic}(0).
$$
Moreover, since $V$ is regular we may assume that
we are in the situation of Example~\ref{ex:canV2}.
By Example~\ref{ex:canV}, sections of
$\EV$ are uniquely written as finite sums 
\begin{equation}
\label{eq:sectEV}
\sum_{\alpha\in\N^d} A_\alpha \partial_{x'}^\alpha, 
\mbox{ with $A_\alpha\in\EX[0]|_V$.}
\end{equation} 
One concludes using \eqref{eq:sectDV} 
since, by definition of $\shl_V$,
$$
\shl_V(\sum_\alpha A_\alpha \partial_{x'}^\alpha) = 
\sum_\alpha\sigma_0(A_\alpha)\partial_{x'}^\alpha.
$$
\end{proof}

\begin{corollary}
\label{co:KO}
Let $V\subset \dot T^*X$ be a conic regular involutive submanifold, and $\stkt$ be a stack of twisted sheaves on $V$. Then there is an equivalence of categories
$$
\catMod(\EV/\EV[-1];\stkt) \simeq 
\catMod(\D^{bic}_V(0);\stkt\stktens\stks_{\Omega_{V/X}^{-1/2}}),
$$
where $\stks_{\Omega_{V/X}^{-1/2}}$ denotes a stack of twisted sheaves such that $\Omega_{V/X}^{-1/2} \in \catMod(\OV;\stks_{\Omega_{V/X}^{-1/2}})$.
\end{corollary}

\section{Statement of the result}

We can now state our main result. 

\begin{theorem}
\label{th:main}
Let $V\subset \dot T^*X$ be a conic involutive submanifold and 
$\Sigma\subset V$ a bicharacteristic leaf. 
Let $\stkt$ be a stack of twisted sheaves on $X$,
and $\shm$ an object of $\catMod(\EX;\stkopb{\pi}\stkt)$ 
globally simple  along $V$. Then 
$$
\pi_\Sigma^\sharp(\stkclass\stkt) = 
\stkclass{\stks_{\Omega^{1/2}_{\Sigma/X}}} \quad
\text{in }H^2(\Sigma;\C_\Sigma^\times),
$$
where $\pi_\Sigma^\sharp \colon H^2(X;\CX^\times) \to H^2(\Sigma;\C_\Sigma^\times)$ denotes
the pull-back and
$\stks_{\Omega_{\Sigma/X}^{1/2}}$ denotes a stack of twisted sheaves such that $\Omega_{\Sigma/X}^{1/2} \in \catMod(\O_\Sigma;\stks_{\Omega_{\Sigma/X}^{1/2}})$.
\end{theorem}

\begin{proof}
The proof follows the same lines as in \cite[\S I.5.2]{KK}.
Let us first reduce to the regular involutive case by the trick of the dummy variable. 
Let $p\colon \widetilde X = X\times\C \to X$ be the projection. 
With the notations of Proposition~\ref{pr:dummy}, 
replace $X$ with $\widetilde X$, $\stkt$ with 
$\widetilde\stkt = \stkopb p \stkt$, $V$ with 
$\widetilde V = V \times \dot T^*\C$,  $\shm$ with 
$\widetilde \shm$, and $\Sigma$ with 
$\widetilde\Sigma = \Sigma \times \{(0;1)\}$. 
Under the isomorphism 
$H^2(\Sigma;\C_\Sigma^\times)\simeq 
H^2(\widetilde\Sigma;\C_{\widetilde\Sigma}^\times)$ one has 
$\pi_\Sigma^\sharp(\stkclass\stkt)=\pi_{\widetilde\Sigma}^\sharp(\stkclass{\widetilde\stkt})$. 
Hence we may assume that $V$ is regular involutive.

Let $\shm_0$ be a $V$-lattice in $\shm$. 
By definition of twisted global simplicity, 
$\overline\shm_0 = \shm_0/\filt[-1]\shm$ is an object of 
$\catMod(\EV/\EV[-1];\stkopb{\pi_V}\stkt)$ locally 
isomorphic to $\OV[0]$. By Theorem~\ref{th:LV1}, 
one has a $\C$-equivalence
$$
\catMod(\EV/\EV[-1];\stkopb{\pi_V}\stkt)
\simeq \catMod(\D^{bic}_{\Omega_{V/X}^{1/2}}(0);\stkopb{\pi_V}\stkt).
$$
Denote by $j_\Sigma\colon \Sigma\to V$ the 
embedding of the bicharacteristic leaf. 
Then $\dzeroopb{j_\Sigma}(\overline\shm_0)$ is an object of 
$\catMod(\D_{\Omega_{\Sigma/X}^{1/2}};\stkopb{\pi_\Sigma}\stkt)$ 
locally isomorphic to $\O_\Sigma$.
Using \eqref{eq:DF-SeqvD-FS}, we get 
the equivalence of $\C$-stacks
\begin{align*}
\catMod(\D_{\Omega_{\Sigma/X}^{1/2}};\stkopb{\pi_\Sigma}\stkt) &\isoto
\catMod(\D_\Sigma;\stkopb{\pi_\Sigma}\stkt\stktens\stks_{\Omega_{\Sigma/X}^{-1/2}}).
\end{align*}
Since $\dzeroopb{j_\Sigma}(\overline\shm_0)$ 
is a flat connection  of rank $1$ in 
$\catMod(\D_\Sigma;\stkopb{\pi_\Sigma}\stkt\stktens\stks_{\Omega_{\Sigma/X}^{-1/2}})$, its solution sheaf
$\hom[\D_\Sigma](\dzeroopb{j_\Sigma}(\overline\shm_0),\O_\Sigma)$
is a local system of rank $1$
in $(\stkopb{\pi_\Sigma}\stkt\stktens\stks_{\Omega_{\Sigma/X}^{-1/2}})(\Sigma)$.
The statement then follows by Proposition~\ref{pr:trivialstk}.
\end{proof}

\begin{remark}
Let us say that a coherent $\EX$-module 
$\shm$ is globally $r$-simple along $V$ if it admits a lattice
$\shm_0$  such that 
$\E_V\shm_0\subset\shm_0$ and
$\shm_0/\filt[-1]\shm$ 
is locally isomorphic to $\OV(0)^r$. 

Theorem \ref{th:main} extends to globally $r$-simple modules as follows. If  
an object $\shm$ of $\catMod(\EX;\stkopb{\pi}\stkt)$ is
globally $r$-simple along $V$, then 
$$
\pi_\Sigma^\sharp(\stkclass{\stkt})^r = 
(\stkclass{\stks_{\Omega^{1/2}_{\Sigma/X}}})^r \quad
\text{in }H^2(\Sigma;\C_\Sigma^\times).
$$
The proof goes along the same lines as the one above, recalling the following fact. Let $\stks$ be a stack of twisted sheaves on $X$, and let $F\in\stks(X)$ be a local system of rank $r$. Then $\det F$ is a local system of rank $1$ in $\stkinv[r]\stks(X)$, so that $\stkinv[r]\stks$ is globally $\C$-equivalent to $\stkMod(\CX)$.

\end{remark}

\begin{corollary}
\label{co:main}
Let $V\subset \dot T^*X$ be a conic involutive submanifold and 
$\Sigma\subset V$ a bicharacteristic leaf.
Let $\stkt$ be a stack of twisted sheaves on $X$
and $\shm$ an object of $\catMod(\EX;\stkopb{\pi}\stkt)$ 
globally simple  along $V$.
Assume that $\pi_\Sigma^\sharp \colon H^2(X;\CX^\times) \to
H^2(\Sigma;\C_\Sigma^\times)$ is injective and 
that $\stkclass{\stks_{\Omega^{1/2}_{\Sigma/X}}} = 1$ in 
$H^2(\Sigma;\C_\Sigma^\times)$. Then $\stkt$ is globally
$\C$-equivalent to $\stkMod(\C_X)$.
\end{corollary}

\begin{proof}
By Theorem \ref{th:main}, 
$\pi_\Sigma^\sharp(\stkclass\stkt)=1$ in $H^2(\Sigma;\C_\Sigma^\times)$.
Since $\pi_\Sigma^\sharp \colon H^2(X;\CX^\times) \to
H^2(\Sigma;\C_\Sigma^\times)$ is injective, $\stkclass\stkt=1$ in 
$H^2(X;\CX^\times)$, and this implies that the stack $\stkt$ is globally
$\C$-equivalent to $\stkMod(\C_X)$.
\end{proof}

\section{Application: non existence of twisted wave equations}

Let $\TT$ be an $(n+1)$-dimensional complex vector space,
$\PP$ the projective space of lines in $\TT$, and
$\GG$ the Grassmannian of $(p+1)$-dimensional
subspaces  in $\TT$. Assume $n\geq 3$ and $1\leq p\leq n-2$.
The Penrose correspondence (see~\cite{EPW}) is associated with the double fibration 
\begin{equation} 
\label{eq:Pen}
\PP\from[f]\FF\to[g]\GG
\end{equation}
where $\FF=\{(y,x)\in\PP\times \GG; y\subset x\}$ 
is the 
incidence relation, and $f$, $g$ are the natural
projections.
The double fibration \eqref{eq:Pen} induces the maps
$$
\dot T^*\PP \from[p] \dot T^*_\FF(\PP\times\GG)\to[q] \dot T^*\GG,
$$
where $T^*_\FF(\PP\times\GG) \subset T^*(\PP\times\GG)$ denotes the conormal bundle to $\FF$, and $p$ and $q$ are the natural projections. 
Note that $p$ is smooth surjective, and $q$ is a closed embedding.
Set 
$$
V = q(\dot T^*_\FF(\PP\times\GG)).
$$
Then $V$ is a closed conic regular involutive submanifold of $\dot T^*\GG$, and $q$ identifies the fibers of $p$ with the bicharacteristic leaves of $V$. 

For $m\in\Z$, let $\OP(m)$ be the line bundle
on $\PP$ corresponding to the sheaf of 
homogeneous functions of degree $m$ on $\TT$, and denote by
$\shn_{(m)} \eqdot\DP\tens[\O]\OP(-m)$ the associated $\DP$-module.
Denote by $\doim g$ and $\dopb f$ the direct and inverse 
image in the derived categories 
of $\D$-modules and consider the family of $\DG$-modules
$$
\shm_{(1+m/2)} \eqdot H^0 (\doim g \dopb f \shn_{(m)} ).
$$
For $n=3$ and $p=1$, Penrose identifies $\GG$ with a conformal compactification of the complexified Minkowski space, and the $\DG$-module $\shm_{(1+m/2)}$ corresponds to the massless field equation of helicity $1+m/2$.

By \cite{DS1}, for $m\in\Z$, the microlocalization $\E_\GG\tens[\opb\pi\D_\GG]\opb\pi\shm_{(1+m/2)}$ of $\shm_{(1+m/2)}$ is globally simple along $V$. 

\begin{theorem}
Let $\stks$ be a stack of twisted sheaves on $\GG$
and $\shm$ an object of $\catMod(\DG;\stks)$ whose microlocalization
$\E_\GG\tens[\opb\pi\D_\GG]\opb\pi\shm$
is globally simple along $V$.
Then $\stks$ is globally $\C$-equivalent to $\stkMod(\C_\GG)$, so that 
$\catMod(\DG;\stks)$ is $\C$-equivalent to $\catMod(\DG)$.
\end{theorem}

In other words, $\shm$ is untwisted.

\begin{proof}
Let us start by recalling the microlocal geometry underlying the double fibration \eqref{eq:Pen}.
There are identifications 
\begin{align*}
T^*\PP &= \{(y;\eta); y\subset \TT,\eta\in\Hom(\TT/y,y)\}, \\
T^*\GG &= \{(x;\xi); x\subset \TT,\xi\in\Hom(\TT/x,x)\}, \\
T^*_\FF(\PP\times\GG) &= \{(y,x;\tau); y\subset x\subset \TT,\tau\in\Hom(\TT/x,y)\}.
\end{align*}
The maps $p$ and $q$ are described as follows:
$$
\xymatrix@R=0pt{
\dot T^*\PP & \dot T^*_\FF(\PP\times\GG) \ar[l]^-p \ar[r]_-q & \dot T^*\GG \\
(y;\tau \circ j) & (y,x;\tau) \ar@{|->}[l] \ar@{|->}[r] & (x; i \circ \tau),
}
$$
where $i\colon y\rightarrowtail x$ and $j\colon \TT/y \twoheadrightarrow \TT/x$ are the natural maps.
We thus get
$$
V=\{(x;\xi); \rk(\xi)=1\},
$$
where $\rk(\xi)$ denotes the rank of the linear map $\xi$.
In order to describe the bicharacteristic leaves of $V$, denote by $\bb$ the dual projective space consisting of hyperplanes $z\subset\TT$, and consider the incidence relation
$$
\AA = \{ (y,z)\in \PP \times \bb; y\subset z\subset\TT \}. 
$$
Then
$$
\dot T^*_\AA(\PP\times\bb) = \{(y,z;\theta); y\subset z\subset \TT,\ \theta\colon \TT/z \isoto y\}.
$$
There is an isomorphism
\begin{align*}
\dot T^*_\AA(\PP\times\bb) & \isoto \dot T^*\PP \\
(y,z;\theta) &\mapsto (y; \theta \circ k).
\end{align*}
where $k\colon \TT/y \twoheadrightarrow \TT/z$ is the natural map. 
Set $y = \im\xi$, $z = x+\ker\xi$, and consider the commutative diagram of linear maps
$$
\xymatrix{
\TT/y \ar[r]^j \ar[dr]^k & \TT/x \ar[d]^-\ell \ar[r]^\xi & y \ar[r]^i & x \\
& \TT/z \ar[ur]_-{\widetilde\xi} .
}
$$
We thus get the following description of the composite map
$$
\xymatrix@R=0pt@C=2ex{
\widetilde p\colon & V \ar[r]^-\sim & \dot T^*_\FF(\PP \times \GG) \ar[r]^-p & \dot T^*\PP \ar[r]^-\sim & \dot T^*_\AA(\PP \times \bb)\\
& (x;\xi) \ar@{|->}[r] & (\im\xi,x;\xi) \ar@{|->}[r] & (\im\xi;\xi \circ j) \ar@{|->}[r] & (\im\xi, x+\ker\xi;\widetilde\xi), 
}
$$
It follows that the bicharacteristic leaf $\Sigma_{(y,z,\theta)} \eqdot \widetilde p^{-1}(y,z,\theta)$ of $V$ is given by
\begin{align}
\label{eq:Sigma}
\Sigma_{(y,z,\theta)} &= \{(x;\xi); y = \im\xi,\  z = x+\ker\xi,\ \theta\circ\ell = \xi \} \\
\nonumber
&= \{(x;\xi); y\subset x\subset z,\ \xi= \theta\circ 
\ell\},
\end{align}
where $\ell\colon \TT/x 
\twoheadrightarrow \TT/z$ is the natural map. 
Thus, $\Sigma_{(y,z,\theta)}$ is the Grassmannian of $p$-dimensional linear subspaces in the $(n-1)$-dimensional vector space $z/y$. 

Let us fix a point $(y,z,\theta) \in \dot T^*_\AA(\PP\times\bb)$, and set $\Sigma=\Sigma_{(y,z,\theta)}$.
In order to apply Corollary \ref{co:main}, we need to compute the map $\pi_\Sigma^\sharp$ and the class $\stkclass{\stks_{\Omega^{1/2}_{\Sigma/\GG}}}$. 

The universal bundle $U_\GG \to \GG$ is the subbundle of the trivial bundle
$\GG\times\TT$ whose fiber at $x\in \GG$ is the 
$(p+1)$-dimensional linear subspace  $x\subset \TT$ itself.
Consider the line bundle $D_\GG = \det U_\GG$, and denote by $\OG(-1)$ the sheaf of its sections. 
Recall the isomorphisms
\begin{align*}
& H^1(\GG;\CG^\times) \simeq H^2(\GG;\OG^\times) \simeq  0 , \\
& H^1(\GG;\OG^\times) \simeq  \Z\mbox{ with generator } \lbclass{\OG(-1)}, \\
& H^1(\GG;\OG^\times/\CG^\times) \simeq H^1(\GG;\OG^\times/\CG^\times) \simeq \C\mbox{ with generator }  \tlbclass{\stkMod(\CG),\OG(-1)}, 
\end{align*}
so that the sequence of abelian groups
$$
H^1(\GG;\CG^\times) \to[\alpha] H^1(\GG;\OG^\times)\to[\beta]
             H^1(\GG;\OG^\times/\CG^\times)\to[\delta] H^2(\GG;\CG^\times)
             \to H^2(\GG;\OG^\times),
$$
is isomorphic to the sequence of additive abelian groups
$$
0 \to \Z \to[\beta] \C \to[\delta] \C/\Z \to 0.
$$
Similar results hold for $\Sigma$, which is also a grassmannian. 
Moreover, $\pi_\Sigma^*\OG(-1)\simeq \O_\Sigma(-1)$
by Lemma~\ref{le:OG} below. Hence $\pi_\Sigma^\sharp$ is the isomorphism
$$
\pi_\Sigma^\sharp \colon H^2(\GG;\C_\GG^\times) \simeq \C/\Z \simeq 
H^2(\Sigma;\C_\Sigma^\times).
$$

There are isomorphisms
$$
\Omega_\GG \simeq \OG(-n-1), \qquad
\Omega_\Sigma \simeq \O_\Sigma(-n+1).
$$
Again by Lemma~\ref{le:OG}, we thus have
$$
\pi_\Sigma^*\Omega_\GG \simeq \pi_\Sigma^*\OG(-n-1) 
\simeq \O_\Sigma(-n-1).
$$
It follows that 
$\Omega_{\Sigma/\GG} \simeq \O_\Sigma(2)$, and thus
$$
\lbclass{\Omega_{\Sigma/\GG}} = 2\quad\text{in } \Z\simeq H^1(\Sigma;\O_{\Sigma}^\times).
$$
Therefore
$$
\tlbclass{\stks_{\Omega^{1/2}_{\Sigma/\GG}},\Omega^{1/2}_{\Sigma/\GG}} = 1\quad\text{in } \C\simeq H^1(\Sigma;\O_{\Sigma}^\times/\C_{\Sigma}^\times),
$$
so that
$$
\stkclass{\stks_{\Omega^{1/2}_{\Sigma/\GG}}} = 
\delta(\tlbclass{\stks_{\Omega^{1/2}_{\Sigma/\GG}},\Omega^{1/2}_{\Sigma/\GG}})= 0 \quad\text{in } \C/\Z\simeq H^2(\Sigma;\C_{\Sigma}^\times).
$$
The statement follows by Corollary~\ref{co:main}.
\end{proof}

\begin{lemma}
\label{le:OG}
There is a natural isomorphism
$\pi_\Sigma^*\OG(-1)\simeq \O_\Sigma(-1)$.
\end{lemma}

\begin{proof}
Recall that $D_\GG$ denotes the determinant of the universal bundle on $\GG$.
Geometrically, we have to prove that there is an isomorphism 
$\delta\colon D_\Sigma \isoto D_\GG|_\Sigma$.

Recall the description \eqref{eq:Sigma}, and let $(x;\xi) \in \Sigma$
for $p=(y,z,\theta)\in \dot T^*\PP$.
Then $(D_\Sigma)_{(x;\xi)} = \det(x/y)$, $(D_\GG)_{(x;\xi)} = \det x$,
and $\delta$ is obtained by a trivialization of $\det y\simeq\C$.
\end{proof}

\providecommand{\bysame}{\leavevmode\hbox to3em{\hrulefill}\thinspace}

\noindent

\vspace*{1cm}

\parbox[t]{15em}
{\scriptsize{
\noindent
Andrea D'Agnolo\\
Universit{\`a} di Padova\\
Dipartimento di Matematica\\
via G. Belzoni, 7;
35131 Padova, Italy\\
dagnolo@math.unipd.it\\
http://www.math.unipd.it/\~{}dagnolo/}}
\qquad
\parbox[t]{15em}
{\scriptsize{
Pierre Schapira\\
Universit{\'e} Pierre et Marie Curie\\
Institut de Math{\'e}matiques\\
175, rue du Chevaleret; 75013 Paris France\\
schapira@math.jussieu.fr\\
http://www.math.jussieu.fr/\~{}schapira/}}

\end{document}